\documentclass[12pt]{amsart}
\usepackage{amssymb}
\usepackage{amsmath}
\usepackage{mathtools}
\usepackage[all]{xy}
\usepackage{enumitem}
\usepackage{xcolor}

 \rm

\renewcommand{\(}{\left(}
\renewcommand{\)}{\right)}

\newcommand{\<}{\langle}
\renewcommand{\>}{\rangle}

\newcommand{\abs}[1]{\left\lvert#1\right\rvert}

\newcommand{\st}{\:|\:}

\newcommand{\C}{{\mathbb{C}}}
\newcommand{\R}{{\mathbb{R}}}

\renewcommand{\phi}{\varphi}

\newcommand{\twisted}{\,\natural \,}

\renewcommand{\H}{{\mathcal{H}}}

\newcommand{\BH}{{\mathcal{B}}(\H)}

\theoremstyle{plain}
\newtheorem{thm}{Theorem}[section]
\newtheorem{lem}[thm]{Lemma}

\theoremstyle{definition}
\newtheorem{defn}[thm]{Definition}

\theoremstyle{remark}

\title{The Weyl Transform of a smooth measure on a real-analytic submanifold}

\author{Mansi Mishra and M.~K.~Vemuri}
\address{Department of Mathematical Sciences, IIT(BHU), Varanasi 221005.}
\email{mansimishra1312@gmail.com}

\begin{document}

\begin{abstract}
If $\mu$ is a smooth measure supported on a real-analytic submanifold of 
$\R^{2n}$ which is not contained in any affine hyperplane, 
then the Weyl transform of $\mu$ is a compact operator.
\end{abstract}

\keywords{Absolute continuity; Compact operator; Real-analytic submanifold;
Twisted convolution.}
\subjclass[2010]{22D10, 22E30, 43A05, 43A80, 53D55}

\maketitle
\thispagestyle{empty}

\section{Introduction}

Let $f \in L^1(\R^n)$.  The Fourier transform of $f$, denoted by $\hat{f}$,
is given by
$$
\hat{f}(\xi) = \int_{\R^n} f(x) e^{-2\pi i x\cdot \xi}\, dx, \quad \xi \in \R^n.
$$
The Riemann-Lebesgue lemma \cite[Theorem 1.2]{weiss} states that if
$f \in L^1(\R^n)$, then
$$
\lim_{\abs{\xi} \rightarrow \infty}\hat{f}(\xi) = 0.
$$
More generally, the Fourier transform of a finite Borel measure $\lambda$ on 
$\R^n$ is given by
$$
\hat{\lambda}(\xi) = \int_{\R^n} e^{-2\pi i x\cdot \xi}\, d\lambda(x),
\quad \xi \in \R^n.
$$
If $\lambda$ is absolutely continuous with respect to the Lebesgue measure $m$
on $\R^n$, i.e., $\lambda=fm$ for some $f \in L^1(\R^n)$, then the formula above
reduces to the usual definition of the Fourier transform of $f$,
and
$$
\lim_{\abs{\xi} \rightarrow \infty}\hat{\lambda}(\xi) = 0.
$$
In general, the Fourier transform of a measure need not vanish at infinity. 
For example, the Fourier transform of $\delta_0$, the Dirac measure, is
identically equal to $1$. 
However, the decay of the Fourier transform of a measure can be deduced from
certain curvature properties of the support of the measure.

Suppose $M$ is a smooth submanifold of $\R^n$. By a {\em smooth measure on $M$}
we mean a measure of the form $\mu=\psi\sigma$ where $\sigma$ is the measure
on $M$ induced by the Lebesgue measure on $\R^n$, and $\psi$ is a smooth
function on $\R^n$ whose support intersects $M$ in a compact set.

It is well known that if $\mu$ is a smooth measure on a hypersurface in $\R^n$,
$n \geq 2$, whose Gaussian curvature is nonzero everywhere, then
\begin{equation*}
\abs{\widehat{\mu}(\xi)} \leq A \abs{\xi}^{(1-n)/2},
\end{equation*}
where $A$ is a constant independent of $\xi$ 
(see e.g. \cite[p348, Theorem 1]{stein} and \cite[Theorem 7.7.14]{hormander}).

In \cite{wtom}, an analogue of this result for the Weyl transform was proved.
Let $\H = L^2\(\R^n\)$, and $\BH$ the set of bounded operators on $\H$.
If $f \in L^1(\R^{2n})$, the {\em Weyl Transform} of $f$ is the operator 
$W(f) \in \BH$ defined by 
\begin{equation*}
(W(f)\phi)(t)
=
\int_{\R^n}\int_{\R^n} f(x,y) e^{\pi i(x \cdot y+2y \cdot t)} \phi(t+x) \, dx \, dy.  
\end{equation*}
More generally, if $\lambda$ is a finite Borel measure on $\R^{2n}$,
the {\em Weyl Transform} of $\lambda$ is the operator $W(\lambda) \in \BH$ 
defined by 
\begin{equation*}
(W(\lambda)\phi)(t)
=
\int_{\R^{2n}} e^{\pi i(x \cdot y+2y \cdot t)} \phi(t+x) \, d\lambda(x,y).  
\end{equation*}

The Weyl transform can be expressed in terms of the canonical representation
of the Heisenberg group.
Recall that the reduced Heisenberg group $G$ is the set of triples
\begin{equation*}
\{(x,y,z) | x,y \in \R^n, z \in \C, \abs{z} = 1\}
\end{equation*}
with multiplication defined by
\begin{equation*}
(x,y,z)(x',y',z')=\(x+x', y+y', zz'e^{\pi i (x \cdot y'-y \cdot x')}\).
\end{equation*}
According to the Stone-von Neumann Theorem \cite[Theorem 1.50]{folland},
there is a unique irreducible unitary representation $\rho$ of $G$ such that
\begin{equation*}
\rho(0,0,z)= zI.
\end{equation*}
The standard realization of this representation is on the Hilbert space $\H$ by
the action
\begin{equation}\label{E:schrodinger}
\(\rho(x,y,z)\phi\)(t)= ze^{\pi i (x\cdot y+2y\cdot t)}\phi(t+x).
\end{equation}
Thus the Weyl transform of $\lambda$ may be expressed as
\begin{equation}\label{E:weyl}
W(\lambda)
=
\int_{\R^{2n}} \rho(x,y,1) \, d\lambda(x,y),  
\end{equation}
where the integral is the weak integral as defined in
\cite[Definition 3.26]{rudin-fa}.

The analogue of the Riemann-Lebesgue lemma for the Weyl transform is the
fact that $W(f)$ is a compact operator if $f \in L^1(\R^{2n})$
(see e.g. \cite[Theorem 1.30]{folland} and 
\cite[Theorem 1.3.3]{thangavelu-hahg}).
If $\lambda$ is absolutely continuous with respect to the Lebesgue measure
$m$ on $\R^{2n}$, i.e., $\lambda=fm$ for some $f \in L^1(\R^{2n})$, then 
$W(\lambda)$ reduces to the usual definition of the Weyl transform of $f$,
and hence $W(\lambda)$ is a compact operator.

The main result of \cite{wtom} is the following theorem.
\begin{thm}\label{T:wtom}
Suppose $S$ is a compact connected smooth hypersurface in $\R^{2n}$, $n \geq 2$,
whose Gaussian curvature is positive everywhere.
Let $\mu$ be a smooth measure on $S$.
Then $W(\mu)$ is a compact operator.  Moreover, when $n\geq 6$ and $p > n$, 
\begin{equation*}
W(\mu) \in S^{p}(\H),
\end{equation*}
where $S^{p}(\H)$ denotes the $p^{\mathrm{th}}$ Schatten class of $\H$.
\end{thm}

Naturally, the question of what happens if we take a submanifold of arbitrary
codimension into consideration arises.

Let $M$ be a smooth $m$-dimensional submanifold of $\R^n$, $1 \leq m \leq n-1$,
and let $\mu$ be a smooth measure on $M$.
It is well known that if $M$ is of finite 
type, i.e., at each point, $M$ has at most
a finite order of contact with any affine hyperplane, then
\begin{equation*}
\abs{\widehat{\mu}(\xi)} \leq A \abs{\xi}^{-1/k},
\end{equation*}
where $k$ is the type of $M$ inside the support of $\psi$ 
(see \cite[p351, Theorem 2]{stein}). 
In particular,
\begin{equation}\label{E:stein2}
\lim_{\abs{\xi} \to \infty} \widehat{\mu}(\xi)=0.
\end{equation}

If we consider $M$ to be a real-analytic submanifold of $\R^n$, 
then the condition of being finite type is equivalent to $M$ not lying in any 
affine hyperplane.

The main result of this paper is the following theorem, which is an analogue of
Equation (\ref{E:stein2}) for the Weyl transform of a smooth measure supported
on a real-analytic submanifold of finite type, and is proved in
Section \ref{S:proof}.

\begin{thm}\label{T:main-thm}
Suppose $M$ is a connected real-analytic submanifold of $\R^{2n}$ which
is not contained in an affine hyperplane. Let $\mu$ be a smooth measure on $M$.
Then $W(\mu)$ is a compact operator.
\end{thm}

In Section \ref{S:twisted}, we define and study the twisted convolution of
finite Borel measures, which is an essential tool required to prove the
main result.
In Section \ref{S:curve}, we prove that if $n=1$, then the conclusion of 
Theorem \ref{T:main-thm} holds for a submanifold of finite type without
the additional assumption of real analyticity.
This also proves Theorem \ref{T:wtom} partially for $n=1$. 

In \cite{wtom}, the quantum translation of an operator was defined, and 
as an application of Theorem \ref{T:wtom}, it was shown that for 
$n \geq 6$, there exists a non-zero compact operator on $L^2(\R^n)$
whose quantum translates are linearly dependent.
By using a similar argument, as an application of Theorem \ref{T:main-thm},
we get the following result.

\begin{thm}
There exists a non-zero compact operator $T$ on $L^2(\R^n)$ and 
distinct elements $(x_1,y_1), \dots, (x_{4n+1},y_{4n+1}) \in \R^{2n}$ such that
$\{(x_1,y_1)\cdot T, \dots, (x_{4n+1},y_{4n+1})\cdot T\}$
is a linearly dependent set, where
$(x_i,y_i)\cdot T = \rho(x_i,y_i,1)T\rho(x_i,y_i,1)^{-1}$
is the quantum translation of $T$ by $(x_i,y_i)$, $i=1,\dots, 4n+1$.
\end{thm}

\section{Twisted convolution}\label{S:twisted}
Recall that if $f,g \in L^1(\R^{2n})$, the twisted convolution of $f$ and $g$,
denoted by $f \twisted g$ is 
\begin{equation*}
f \twisted g (x,y) = 
\int_{\R^n}\int_{\R^n} f(x-x',y-y')g(x',y')e^{\pi i (x\cdot y'-y\cdot x')}\, dx' \, dy',
\quad (x,y) \in \R^{2n}.
\end{equation*}
Twisted convolution turns $L^1(\R^{2n})$ into a non-commutative Banach
algebra. 
It is well known that the Weyl transform is an algebra homomorphism from
$L^1(\R^{2n})$ to $\BH$, i.e., $W(f\twisted g) = W(f)W(g)$
(see e.g. \cite[p26]{folland} and \cite[p16]{thangavelu-hahg}).

\begin{defn}\label{D:twisted}
Let $\mu$ and $\nu$ be finite Borel measures on $\R^{2n}$. The twisted 
convolution of $\mu$ and $\nu$ is the measure on $\R^{2n}$, denoted by 
$\mu \twisted \nu$, given by
\begin{equation*}
\mu \twisted \nu(E)= 
\int_{\R^{2n}}\int_{\R^{2n}}\chi_{E}(x+x',y+y')e^{\pi i (x\cdot y'-y\cdot x')}\, 
d\mu(x,y) \, d\nu(x',y'),
\end{equation*}
where $E$ is a Borel subset of $\R^{2n}$.
\end{defn}

It follows from Definition \ref{D:twisted} that if $f \in C_0(\R^{2n})$, then
\begin{equation}\label{E:twisted}
\int_{\R^{2n}} f(x,y)\, d(\mu \twisted \nu)(x,y)
=
\int_{\R^{2n}}\int_{\R^{2n}}f(x+x',y+y')e^{\pi i (xy'-yx')}\, 
d\mu(x,y) \, d\nu(x',y').
\end{equation}

Let $f,g \in L^1(\R^{2n})$. Let $m$ denote the Lebesgue measure on $\R^{2n}$.
Let $\mu_f= fm$ and $\mu_g = gm$. 
If $E$ is a Borel subset of $\R^{2n}$, then
\begin{equation*}
\begin{aligned}
\mu_f \twisted \mu_g(E)
&\;= 
\int_{\R^{2n}}\int_{\R^{2n}}\chi_{E}(x+x',y+y')e^{\pi i (x\cdot y'-y\cdot x')}\, 
d\mu(x,y) \, d\nu(x',y')\\
&\;= 
\int_{\R^{2n}}\int_{\R^{2n}}\chi_{E}(x+x',y+y')e^{\pi i (x\cdot y'-y\cdot x')}f(x,y)g(x',y')\, 
dm(x,y) \, dm(x',y')\\
&\;= 
\int_{\R^{2n}}\chi_{E}(x,y)\int_{\R^{2n}}
e^{\pi i (x\cdot y'-y\cdot x')}f(x-x',y-y')g(x',y') \, dm(x',y')\, dm(x,y)\\
&\;= 
\int_{\R^{2n}}\chi_{E}(x,y) (f\twisted g)(x,y)\, dm(x,y)\\
&\;= 
\mu_{f\twisted g}(E).
\end{aligned}
\end{equation*}
Therefore the two definitions of twisted convolution coincide for
$L^1$ functions.

Twisted convolution turns $M(\R^{2n})$, the set of finite Borel measures on
$\R^{2n}$, into a non-commutative Banach algebra.
The following theorem shows that the Weyl transform is an algebra homomorphism
from $M(\R^{2n})$ to $\BH$.

\begin{thm}\label{T:weyl}
Let $\mu$ and $\nu$ be finite Borel measures on $\R^{2n}$.  Then
$$
W(\mu \twisted \nu) = W(\mu)W(\nu).
$$
\end{thm}

\begin{proof}
By Equations (\ref{E:weyl}) and (\ref{E:twisted}),
\begin{equation*}
\begin{aligned}
W(\mu \twisted \nu)
=&\,
\int_{\R^{2n}} \rho(x,y,1) \,d(\mu \twisted \nu)(x,y)\\
=&\,
\int_{\R^{2n}}\int_{\R^{2n}} \rho(x+x',y+y',1) e^{\pi i (x\cdot y' - y\cdot x')}
\,d\mu(x,y)\,d\nu(x',y')\\
=&\,
\int_{\R^{2n}}\int_{\R^{2n}} \rho(x+x',y+y',e^{\pi i (x\cdot y' - y\cdot x')})
\,d\mu(x,y)\,d\nu(x',y')
\quad \text{(by Equation (\ref{E:schrodinger}))}\\
=&\,
\int_{\R^{2n}}\int_{\R^{2n}} \rho((x,y,1)(x',y',1)) \,d\mu(x,y)\,d\nu(x',y')\\
=&\,
\int_{\R^{2n}}\int_{\R^{2n}} \rho(x,y,1)\rho(x',y',1) \,d\nu(x',y') \,d\mu(x,y)\\
=&\,
\int_{\R^{2n}} \rho(x,y,1)W(\nu) \,d\mu(x,y)\\
=&\,
W(\mu)W(\nu).
\end{aligned}
\end{equation*}
\end{proof}

\section{The proof}\label{S:proof}

To prove Theorem \ref{T:main-thm}, we need a result analogues to a result of
Ragozin about the absolute continuity of the convolution of measures supported 
on analytic submanifolds of $\R^n$ (see \cite[Theorem 5.1]{ragozin}).
In particular, Ragozin proved the absolute continuity of the convolution 
square of the surface measure on a compact analytic hypersurface of $\R^n$.
Later, Thangavelu proved that the twisted convolution of the surface measure
on a unit sphere in $\R^{2n}$ with itself is absolutely continuous with
respect to the Lebesgue measure on $\R^{2n}$
(see \cite[Proposition 4.3]{thangavelu}).

The following theorem is an analogue of \cite[Theorem 5.1]{ragozin}
for twisted convolution. 

\begin{thm}\label{T:ragozin}
Let $M_1, \dots, M_k$ be connected real-analytic submanifolds of $\R^{2n}$
such that for some choice of points $p_i \in M_i$, $i=1, \dots, k$,
$T_{p_1}M_1 + \cdots + T_{p_k}M_k = \R^{2n}$, where $T_{p_i}M_i$ denotes the 
tangent space at $p_i$ of $M_i$. If $\mu_i$ is a smooth measure on 
$M_i$, $i=1, \dots, k$, then $\mu_1 \twisted \cdots \twisted \mu_k$ is
absolutely continuous with respect to the Lebesgue measure on $\R^{2n}$.
\end{thm}

Let $\Sigma_k: M_1 \times \dots \times M_k \to \R^{2n}$ be the map given by
$\Sigma_k (p_1,\dots, p_k)= p_1+ \cdots+ p_k$.
We need the following lemma to prove Theorem \ref{T:ragozin}.

\begin{lem}\label{L:pf}
There exists a smooth
function $\phi_k: M_1 \times \dots \times M_k \to \C$ such that,
for a Borel set $E \subseteq \R^{2n}$,
\begin{equation*}
\mu_1 \twisted \cdots \twisted \mu_k (E)= 
\(\phi_k \,\mu_1 \times \dots \times \mu_k\)(\Sigma_k^{-1}(E)),
\end{equation*}
i.e., $\mu_1 \twisted \cdots \twisted \mu_k$ is the push-forward of the measure
$\phi_k \,\mu_1 \times \dots \times \mu_k$ by $\Sigma_k$.
\end{lem}

\begin{proof}
We will prove the result by induction on $k$.

Let $k=2$. Observe that by Definition \ref{D:twisted}, 
\begin{equation}\label{E:2d-twisted}
\mu_1 \twisted \mu_2 (E)= 
\(\phi_2 \,\mu_1 \times \mu_2\)(\Sigma_2^{-1}(E)),
\end{equation}
where $\phi_2: \R^{2n}\times \R^{2n} \to \C$ is given by
$\phi_2((x_1,y_1), (x_2,y_2)) = e^{\pi i \(x_1\cdot y_2 - y_1 \cdot x_2\)}$.

Assume that there exists a function
$\phi_{k-1}: M_1 \times \dots \times M_{k-1} \to \C$ such that 
$$
\mu_1 \twisted \cdots \twisted \mu_{k-1} (E)= 
\(\phi_{k-1} \,\mu_1 \times \dots \times \mu_{k-1}\)(\Sigma_{k-1}^{-1}(E)).
$$

Then
\begin{equation*}
\begin{gathered}
((\mu_1 \twisted \cdots \twisted \mu_{k-1}) \twisted \mu_k) (E)\\
=
\int_{\R^{2n}}\int_{\R^{2n}}\chi_{E}(x+x_k,y+y_k)e^{\pi i (x\cdot y_k-y\cdot x_k)}\, 
d(\mu_1 \twisted \cdots \twisted \mu_{k-1})(x,y) \, d\mu_k(x_k,y_k)\\
=
\int_{\R^{2n}}\dots\int_{\R^{2n}}
\left[ \chi_{E}(x_1+ \dots +x_k,y_1 +\dots +y_k)
e^{\pi i ((x_1+\cdots +x_{k-1})\cdot y_k-(y_1+\cdots +y_{k-1})\cdot x_k)} \right.\\
\left.\phi_{k-1} ((x_1,y_1),\dots, (x_{k-1},y_{k-1})) \right]
\, d\mu_1(x_1,y_1) \dots \, d\mu_k(x_k,y_k)\\
=
(\phi_k \,\mu_1 \times \dots \times \mu_k)(\Sigma_k^{-1}(E)),
\end{gathered}
\end{equation*}
where $\phi_k: M_1 \times \dots \times M_k \to \C$ is given by
$$
\phi_k((x_1,y_1), \dots, (x_k,y_k))= 
\phi_{k-1}((x_1,y_1), \dots, (x_{k-1},y_{k-1}))
e^{\pi i ((x_1+\cdots +x_{k-1})\cdot y_k-(y_1+\cdots +y_{k-1})\cdot x_k)}.
$$
\end{proof}

Theorem \ref{T:ragozin} now follows by the argument in 
\cite[Theorem 5.1]{ragozin}.  Here, we give a more streamlined proof using the 
coarea formula. 

Let $M_1, \dots, M_k$ be connected real-analytic submanifolds of $\R^{2n}$.
Then $M_1 \times \dots \times M_k$ is a connected real-analytic manifold.
Let $\tau_{(M_1 \times \dots \times M_k)}$ denote the Riemannian measure on
$M_1 \times \dots \times M_k$.

Observe that $\Sigma_k$ is an analytic map.
Let $p_i \in M_i$, $i=1, \dots, k$ be such that
$T_{p_1}M_1 + \cdots + T_{p_k}M_k = \R^{2n}$.
Then the rank of
$\Sigma_k$ is $2n$ at the point $(p_1, \dots, p_k)$.
Therefore the critical set of
$\Sigma_k$ is a proper analytic subvariety of $M_1 \times \dots \times M_k$,
and hence has $\tau_{(M_1 \times \dots \times M_k)}$ measure zero.

If $\mu_i$ is a smooth measure on $M_i$, $i=1, \dots, k$, then
$\mu_1 \times \dots \times \mu_k$ is a smooth measure on 
$M_1\times \dots \times M_k$, and so
$\phi_k \, \mu_1 \times \dots \times \mu_k$ is a smooth measure on
$M_1\times \dots \times M_k$.

The proof of Theorem \ref{T:ragozin} now follows from the following lemma.

\begin{lem}\label{L:push-forward}
Let $M$ and $N$ be Riemannian manifolds. Let $\tau$ and $\nu$ denote the
Riemannian measures on $M$ and $N$ respectively.
Let $\mu= \psi \tau$ be a smooth measure on $M$.
Suppose $f:M \to N$ is a differentiable map.
If $f$ is a submersion everywhere except on a set of $\tau$-measure zero,
then the
push-forward of $\mu$ by $f$ is absolutely continuous with respect to $\nu$.
\end{lem}

\begin{proof}
Let $Z=\{x\in M \st \text{$f$ is not a submersion at $x$}\}$.
Then $\tau(Z)=0$, and so $\mu(Z)=0$.
Let $\mathcal{J}_f$ denote the normal Jacobian of $f$, i.e.,
the absolute value of the determinant of $df$ restricted to the
orthogonal complement of its kernel.  Then $\mathcal{J}_f$ is
strictly positive on the set of regular points of $f$, i.e., on $M \setminus Z$.

Let $f_*\mu$ denote the push-forward of $\mu$ by $f$.
For $x \in N$, let $\sigma_x$ denote the Riemannian measure on the manifold
$f^{-1}(x)$.
Let $U \subseteq N$ be a Borel set. 
By the coarea formula (see \cite[p159]{chavel}), we have
\begin{equation*}
\begin{aligned}
f_*\mu (U)
=&\;
\mu(f^{-1}(U))\\
=&\;
\mu(f^{-1}(U)\setminus Z)\\
=&\;
\int_{f^{-1}(U) \setminus Z} d\mu\\
=&\;
\int_{f^{-1}(U) \setminus Z} \psi \, d\tau\\
=&\;
\int_{U \setminus f(Z)} \int_{f^{-1}(x)} 
\frac{\psi(y)}{\mathcal{J}_f(y)} \,d\sigma_x(y) \,d\nu(x),
\end{aligned}
\end{equation*} 
i.e.,
$\frac{df_*\mu}{d\nu}(x)=\int_{f^{-1}(x)} 
\frac{\psi(y)}{\mathcal{J}_f(y)}\,d\sigma_x(y)$,
and so $f_*\mu$ is absolutely continuous with respect to $\nu$.
\end{proof}

We need the following lemma to prove Theorem \ref{T:main-thm}.

\begin{lem}\label{L:ind}
Let $M$ be a connected real-analytic submanifold of $\R^{2n}$.
Assume that $M$ is not contained in any affine hyperplane.
Then there exist points $p_1, \dots, p_k \in M$ such that 
$$
T_{p_1}M + \cdots + T_{p_k}M = \R^{2n}.
$$
\end{lem}

\begin{proof}
Suppose $\sum_{p \in M}T_pM \neq \R^{2n}$.  Then there exists a hyperplane $Y$ in
$\R^{2n}$ such that $\sum_{p \in M}T_pM \subseteq Y$.

Fix $p \in M$. Let $q \in M$. Since $M$ is connected, there exists a
differentiable map $\xi: [0,1] \to M$ such that $\xi(0)=p$ and $\xi(1)=q$.
By the fundamental theorem of calculus,
$$
q=p+ \int_0^1 \dot{\xi}(t)\, dt.
$$
Since $\dot{\xi}(t) \in Y$ for $t \in [0,1]$ and $Y$ is closed,
it follows that $q \in p+ Y$. This implies that $M$ is contained in the affine
hyperplane $p+ Y$, which is a contradiction.

Hence $\sum_{p \in M}T_pM =\R^{2n}$.
By the Steinitz exchange lemma, there exist finitely many points
$p_1, \dots, p_k \in M$ such that $T_{p_1}M +\cdots +T_{p_k}M = \R^{2n}$. 
\end{proof}

By possibly adding one more point, we can assume that the integer $k$ obtained
in Lemma \ref{L:ind} is even.

We are now in a position to prove Theorem \ref{T:main-thm}.
Let $\widetilde{M}= \{-p \st p\in M\}$. 
Then $\widetilde{M}$ is a connected real-analytic submanifold of $\R^{2n}$.
Observe that if $p \in M$, $T_pM = T_{(-p)}\widetilde{M}$.
Therefore, by Lemma \ref{L:ind},
\begin{equation}\label{E:tangentspaces}
T_{p_1}M + T_{(-p_2)}\widetilde{M} +  \cdots + 
T_{p_{k-1}}M + T_{(-p_k)}\widetilde{M} = \R^{2n}.
\end{equation}

Let $\mu$ be a smooth measure on $M$. 
Let $\widetilde{\mu}$ denote the push-forward of $\mu$ by the map
which sends $p\in M$ to $-p$.
Then, by Theorem \ref{T:ragozin} and Equation (\ref{E:tangentspaces}),
it follows that $(\mu \twisted \widetilde{\mu})^{k/2}$ is absolutely continuous 
with respect to the Lebesgue measure on $\R^{2n}$.
Therefore $W((\mu \twisted \widetilde{\mu})^{k/2})$ is a compact operator.

Observe that $W(\widetilde{\mu})= W(\mu)^*$, where $W(\mu)^*$ is the adjoint of 
the operator $W(\mu)$.
By Theorem \ref{T:weyl}, 
$$
W((\mu \twisted \widetilde{\mu})^{k/2}) = (W(\mu)W(\widetilde{\mu}))^{k/2}
= (W(\mu)W(\mu)^*)^{k/2}.
$$
Since $W(\mu)W(\mu)^*$ is self-adjoint and $(W(\mu)W(\mu)^*)^{k/2}$
is compact, it follows that $W(\mu)W(\mu)^*$ is compact.
By the polar decomposition, it follows that $W(\mu)$ is compact.

\section{Curve in $\R^2$}\label{S:curve}
 
In this section, we prove the following result, which states that if $n=1$, the
assumption of real-analyticity can be removed from Theorem \ref{T:main-thm}. 

\begin{thm}\label{T:curve}
Suppose $M \subseteq \R^2$ is a connected smooth hypersurface of finite type. 
Let $\mu$ be a smooth measure on $M$.  Then $W(\mu)$ is a compact operator.
\end{thm}

We need the following lemma to prove this theorem.

\begin{lem}\label{L:2d}
Let $\gamma:[a,b] \to \R^2$ be a finite type unit-speed simple smooth curve.
Let $\delta:[c,d] \to \R^2$ be a unit-speed simple smooth curve.
Let $\mu$ and $\nu$ be smooth measures on $\operatorname{Im}(\gamma)$ and
$\operatorname{Im}(\delta)$ respectively.
Then $\mu \twisted \nu$ is absolutely
continuous with respect to the Lebesgue measure on $\R^2$.
\end{lem}

\begin{proof}
Let $E \subseteq \R^2$ be a Borel set. Then, by Equation (\ref{E:2d-twisted}),
\begin{equation*}
\mu \twisted \nu (E)
=
\(\phi_2 \,\mu \times \nu\)(\Sigma_2^{-1}(E)),
\end{equation*}
where 
$\Sigma_2:\operatorname{Im}(\gamma) \times \operatorname{Im}(\delta) \to \R^2$
is the map given by $\Sigma_2(\gamma(s),\delta(t))= \gamma(s)+\delta(t)$.
Define $S:[a,b]\times [c,d] \to \R^2$ by
$$
S(s,t)= \gamma(s)+\delta(t).
$$
Let 
$T=\{(s,t)\in [a,b]\times [c,d]\st \dot{\gamma}(s)=\pm\dot{\delta}(t)\}$.
Observe that $T$ is the critical set of $S$. 
We claim that $T$ has area zero.
Suppose $T$ has a positive area.
Then, by Fubini's theorem, there exists $t \in [c,d]$ such that the set
$Z=\{s\in [a,b]\st \dot{\gamma}(s)=\pm\dot{\delta}(t)\}$ has positive
length.
By the Lebesgue differentiation theorem \cite[Theorem 3.21]{folland-real},
there exists $s' \in [a,b]$
such that $s'$ is a point of density of $Z$.
Then there exists distinct $s_j \in Z$
such that the sequence $\{s_j\}$ converges to $s'$, and so
the sequence $\{\dot{\gamma}(s_j)\}$ converges to $\dot{\gamma}(s')$.
Let $\vec{v}$ be a unit vector in $\R^2$ perpendicular to $\dot{\delta}(t)$.
Then
$$
\<\dot{\gamma}(s_j)\,,\,\vec{v}\>
=
\<\pm\dot{\delta}(t)\,,\,\vec{v}\>
=
0, \quad j=1,2,\dots.
$$
Hence, all the coefficients in the Taylor expansion of 
$\<\dot{\gamma}(s)\,,\,\vec{v}\>$ about $s=s'$ are zero.
This contradicts the fact that $\gamma$ is a curve of finite type.  
Therefore $T$ has area zero.

Since $\operatorname{Im}(\gamma)$ and $\operatorname{Im}(\delta)$ are smooth
submanifolds of $\R^2$, it follows that 
$\operatorname{Im}(\gamma) \times \operatorname{Im}(\delta)$
is a smooth submanifold of $\R^4$. Let $\tau$ denote the Riemannian measure
on $\operatorname{Im}(\gamma) \times \operatorname{Im}(\delta)$.
Since $\gamma \times \delta$ is a smooth map and the set $T$ has area zero,
it follows that $(\gamma \times \delta)(T)$ has $\tau$-measure zero.

Observe that $(\gamma\times \delta)(T)$ is the critical set of $\Sigma_2$.
The proof now follows from Lemma \ref{L:push-forward}.
\end{proof}

We now prove Theorem \ref{T:curve}.
Suppose $M \subseteq \R^2$ is a hypersurface of finite type. 
Let $\mu$ be a smooth measure on $M$. Let $\widetilde{\mu}$ denote
the push-forward of $\mu$ by the map which sends
$p\in M$ to $-p$. 
It follows from Lemma \ref{L:2d} that $\mu \twisted \widetilde{\mu}$ is 
absolutely continuous with respect to the Lebesgue measure on $\R^2$.
Therefore $W(\mu \twisted \widetilde{\mu})$ is compact.
Observe that $W(\widetilde{\mu})= W(\mu)^*$.
By Theorem \ref{T:weyl}, 
$$
W(\mu \twisted \widetilde{\mu})= W(\mu)W(\widetilde{\mu})= W(\mu)W(\mu)^*.
$$
Since $W(\mu)W(\mu)^*$ is compact, it follows
by the polar decomposition that $W(\mu)$ is compact.
\bibliographystyle{amsplain}
\bibliography{v5-wtomas}

\end{document}